\documentclass[a4paper,12pt,oneside]{amsart}
\usepackage[left=2.7cm,right=2.7cm,top=3.5cm,bottom=3cm]{geometry}
\usepackage[all]{xy}
\usepackage{amsfonts}
\usepackage{amssymb}
\usepackage{amsmath}
\usepackage{latexsym}
\usepackage{amscd}
\input cyracc.def
\font\tencyr=wncyr10
\def\cyr{\tencyr\cyracc}
\newcommand{\ts}{\mbox{\cyr Sh}}
\newcommand{\ov}{\overline}
\newcommand{\Q}{\mathbb{Q}}
\newcommand{\Z}{\mathbb{Z}}

\newcommand{\F}{\mathcal{F}}
\newcommand{\X}{\mathcal{X}}
\newcommand{\ol}{\mathcal{O}}
\renewcommand{\S}{\mathcal{S}}
\newcommand{\Ql}{{\mathbb{Q}_l}}
\newcommand{\Zl}{{\mathbb{Z}_l}}
\newcommand{\Zp}{{\mathbb{Z}_p}}

\newcommand{\Zr}{{\mathbb{Z}_r}}
\newcommand{\Qr}{{\mathbb{Q}_r}}
\newcommand{\N}{\mathbb{N}}

\newcommand{\dl}[1]{\lim_{\buildrel \longrightarrow\over{#1}}}
\newcommand{\sri}{\twoheadrightarrow}
\newcommand{\iri}{\hookrightarrow}
\newcommand{\ri}{\rightarrow}

\renewcommand{\L}{\Lambda}

\newcommand{\wh}{\widehat}

\newcommand{\g}{\gamma}
\newcommand{\G}{\Gamma}

\newcommand{\bF}{\mathbb{F}}
\renewcommand{\k}{\kappa}
\newtheorem{thm}{Theorem}[section]
\newtheorem{defin}[thm]{Definition}
\newtheorem{prop}[thm]{Proposition}
\newtheorem{ex}[thm]{Example}

\newtheorem{lemma}[thm]{Lemma}
\newtheorem{cor}[thm]{Corollary}
\newtheorem{rem}[thm]{Remark}
\newcommand{\liminv}{\displaystyle \lim_{\leftarrow}}
\newcommand{\limdir}{\displaystyle \lim_{\rightarrow}}

\begin{document}

\title[Selmer groups ...]{Selmer groups for elliptic curves in $\Z_l^d$-extensions of function fields of
characteristic $p$}

\author{A. Bandini, I. Longhi}

\begin{abstract} Let $F$ be a function field of characteristic $p>0$,
$\F/F$ a $\Z_l^d$-extension (for some prime $l\neq p$) and $E/F$ a
non-isotrivial elliptic curve. We study the behaviour of the $r$-parts of the Selmer groups 
($r$ any prime) in the subextensions of $\F$ via appropriate
versions of Mazur's Control Theorem. As a consequence we prove that the limit of the Selmer groups 
is a cofinitely generated (in some cases cotorsion) module over the Iwasawa algebra of $\F/F$.

\vspace{.5truecm}

\noindent {\sc R\'esum\'e.} Soit $F$ un corps de fonctions de caract\'eristique $p>0$, $\mathcal F/F$
une $\mathbb Z_l^d$-extension (pour un nombre premier $l\neq p$) et $E/F$ une courbe elliptique non-isotrivale.
Nous \'etudions le comportement des $r$-parties des groupes de Selmer pour
les sous-extensions de $\mathcal F$ par des variantes du Th\'eor\`eme de
contr\^ole de Mazur. Cons\'equemment, nous d\'emontrons que la limite des
groupes de Selmer est un module finiment coengendr\'e (parfois de cotorsion) 
sur l'alg\`ebre d'Iwasawa de $\mathcal F/F$.
\end{abstract}

\maketitle

Groupes de Selmer, courbes elliptiques, corps de fonctions, th\'eorie d'Iwasawa.


\section{Introduction}
Let $F$ be a function field (in the whole paper function field means a field
of trascendence degree 1 over its constant field) with constant field
$\bF$ an intermediate extension between $\bF_p$ (the field with $p$ elements)
and a (fixed) algebraic closure $\ov{\bF}_p$ of $\bF_p\,$.
Let $E/F$ be a non-isotrivial elliptic curve (i.e., $j(E)\notin \bF$)
and assume that $E$ has good or split multiplicative reduction
at all primes of $F$ (it is always possible to reduce to this
situation by simply taking a finite extension of $F$). \\
Let $l$ be a prime different from $p$, let $\F/F$ be a $\Z_l^d$-extension of $F$
with Galois group $\G$ (the case $l=p$ has been developed in \cite{BL} for global
function fields). Denote by $\L:=\Zl[[\G]]$ the associated Iwasawa algebra.
Let $\bF_p^{(l)}$ be the unique $\Z_l$-extension of $\bF_p\,$. If $\bF_p^{(l)}\not\subset \bF$
then there is only one $\Z_l$-extension of $F$, namely the arithmetic one, obtained by extending scalars
from $\bF$ to $\bF_p^{(l)}\bF$ (see Proposition \ref{ZlExt}); we recall that this extension is
everywhere unramified. On the other hand, if, for example, $\bF$ contains $\boldsymbol\mu_{l^\infty}$
(the roots of unity of $l$-power order) then Kummer theory produces lots of examples of
disjoint $\Z_l$-extensions of $F$ (see the Appendix).

In section 2 we will define the $r$-part ($r$ any prime) of the Selmer group of $E$,
$Sel_E(L)_r\,$, for any algebraic extension $L$ of $F$. Our goal
is to study the structure of $Sel_E(\F)_r$ (actually of its Pontrjagin
dual) as a $\Z_r[[\G]]$-module.

Not surprisingly the most interesting case happens to be $r=l$.
Let $\S$ be the Pontrjagin dual of $Sel_E(\F)_l\,$: its structure depends, among other things, on the
base field $\bF$. Namely we have different results depending on whether $\bF_p^{(l)}\subset \bF$ or not.
In section \ref{SecControll}, we shall prove the following

\begin{thm}\label{IntrThm1}
Assume that $\bF$ does not contain $\bF_p^{(l)}\,$. Then $\S$ is a finitely generated $\L$-module. Moreover if
$Sel_E(F)_l$ is finite then $\S$ is $\L$-torsion.
\end{thm}

\begin{thm}\label{IntrThm2}
Assume that only finitely many primes of $F$ are ramified in $\F/F$ and
that $\bF$ contains $\bF_p^{(l)}\,$.
Then $\S$ is a finitely generated $\L$-module.\\
Moreover if:\begin{itemize}
\item[1.] the ramified primes are of good reduction for $E$;
\item[2.] for any ramified prime $v$, $E[l^\infty](F_v)$ is finite ($F_v$ is the completion of $F$ at $v$);
\item[3.] $Sel_E(F)_l$ is finite,
\end{itemize}
then $\S$ is $\L$-torsion.
\end{thm}

\begin{rem}\label{IntrRem}
{\em When $F$ is a global function field, according to the Birch and Swinnerton-Dyer conjecture,
$Sel_E(F)_l$ is finite if and only if ${\rm rank}\,E(F)=0$.}
\end{rem}

When $\F/F$ is a $\Z_l$-extension and $\S$ is $\L$-torsion is quite easy to prove that $E(\F)$ is finitely
generated (see Corollary \ref{FinGenE}). The behaviour of the rank of $E$ in
an infinite tower of extensions of a function field $K$ (in any characteristic) has been addressed by many authors.
Among others, Shioda \cite{Sh0}, Fastenberg \cite{Fa} and Silverman \cite{Si3} have provided examples of elliptic curves
with bounded rank in towers of function fields in characteristic 0 and Ulmer \cite{Ul} gives instances of the same phenomenon
for elliptic curves over $\ov{\bF}_q(t^{1/r^m})$ ($r$ a prime not dividing $q$). In the opposite direction examples of
elliptic curves with unbounded rank have been given by Shioda \cite{Sh0} for the tower $\ov{\bF}_p(t^{1/r^m})$
and Ulmer \cite{Ul0} for $\bF_p(t^{1/r^m})$.
In the same spirit the structure of Selmer groups has been studied
by Ellenberg \cite{El} from a slightly different (more geometric) viewpoint
using formulas on Euler characteristic for $\L$-modules.

Since Mazur's classical work \cite{Ma}, duals of Selmer groups have provided the algebraic counterpart for $p$-adic $L$-functions
in Iwasawa theory of elliptic curves over number fields. In section \ref{IMC} we speculate about such an application of our results
when $F$ is a global field.\\

The main tools for the proofs of Theorems \ref{IntrThm1} and \ref{IntrThm2} are appropriate versions of Mazur's
Control Theorem (originally proved in \cite{Ma}; for a different approach, closer to ours, see \cite{Gr1}
and \cite{Gr2}), which we prove in section \ref{SecControll} as well, and Theorem
\ref{NakBH}, a generalization of Nakayama's Lemma which has been proved in \cite{BH}.
We follow some of the basic ideas developed in \cite{BL} for the case $l=p$.

Moreover we can prove a version of the control theorem for
$Sel_E(\F)_r$ for $r\neq l$ as well, but, unfortunately, $Sel_E(\F)_r$ is a module over $\Z_r[[\G]]$,
a ring which we know very little about. Nevertheless we can say something on
the structure of $Sel_E(\F)_r$ and we gathered the results on that module
in section \ref{SecControlp}.

The paper ends with a short Appendix which provides a classification of $\Z_l^d$-extensions of
a field $F$ containing $\boldsymbol\mu_{l^\infty}\,$.\\

\noindent{\bf Acknowledgements.} The authors would like to thank S. Petersen for comments on an
earlier version and for pointing out the idea for the proof of Proposition \ref{ZlExt},
F. Trihan and F. Andreatta for helpful suggestions and discussions.
We are grateful to the anonymous referee for comments which led to improvements in the paper.

\section{The setting and the Selmer groups}
\subsection{Notations} We list some notations which will be used throughout the paper and
briefly describe the setting in which the theory will be developed.

\subsubsection{Fields}\label{FieldsNotations} Let $L$ be a field: then $L^{sep}$ will denote a
separable algebraic closure of $L$ and we put $G_L:=Gal(L^{sep}/L)$.
Moreover $\ov{L}$ will denote an algebraic closure of $L$.\\
If $L$ is a global field (or an algebraic extension of such),
$\mathcal M_L$ will be its set of places. For any place
$v\in\mathcal M_L$ we let $L_v$ be the completion of $L$ at $v$,
$\ol_v$ the ring of integers of $L_v\,$,  $ord_v$ the valuation associated
to $v$ and $\mathbb L_v$ the residue field.\\
As usual, $\boldsymbol\mu_n$ denotes the group of $n$-th roots of 1.

As stated in the introduction, we fix a function field $F$ of
characteristic $p>0$ and an algebraic closure $\ov{F}$. Its constant field will be denoted by
$\bF$. Then $F$ is generated over $\bF$ by a finite number of trascendental elements $z_0\,,\dots ,
z_n$ subjected to algebraic relations. These relations are defined over some finite field
$\bF_q\subset \bF$ for $q\gg 0$. Let $F_0:=\bF_q(z_0\,,\dots ,z_n)$: then $F_0$ is a global field,
$F=\bF F_0$ and $Gal(F/F_0)\simeq Gal(\bF/\bF_q)$.\\
For any place $v\in\mathcal M_F$
we choose $\ov{F_v}$ and an embedding $\ov{F}\hookrightarrow\ov{F_v}$, so to get a corresponding inclusion
$G_{F_v}\hookrightarrow G_F\,$. All algebraic extensions of $F$ (resp. of $F_v$) will be assumed to be contained
in $\ov{F}$ (resp. in $\ov{F_v}$).\\
Script letters will denote infinite extensions of $F$; in particular
$\F/F$ will be a $\Z_l^d$-extension with $l$ a fixed prime different
from $p$. We shall consider a sequence of finite extensions of $F$
such that
\[ F\subset F_1\subset\dots\subset F_n\subset\dots\subset
\bigcup F_n=\F \ .\] In this setting we let $\G:=Gal(\F/F)$ and
$\G_n:=Gal(\F/F_n)$ (for any $n > 0$).

For $\g$ an element in a profinite group, $\ov{\langle\g\rangle}$ will denote
the closed subgroup topologically generated by $\g$.

\subsubsection{Elliptic curves} We fix a non-isotrivial elliptic curve $E/F$,
having split multiplicative reduction at all
places supporting its conductor. The reader is reminded that then at such
places $E$ is isomorphic to a Tate curve, i.e. $E(F_v)\simeq F_v^{\ast}/q_{E,v}^{\Z}$ for
some $q_{E,v}$ (the \emph{Tate period} at $v$) with $ord_v(q_{E,v})=-ord_v(j(E))>0$.\\
For any positive integer $n$ let
$E[n]$ be the scheme of $n$-torsion points. Moreover, for any prime $r$, let
$E[r^\infty]:=\limdir E[r^n]\,$.\\
By the theory of the Tate curve, if $v$ is of bad reduction for $E$ and $r\neq p$ one has
an isomorphism of Galois modules
\[ E[r^\infty](\ov{F_v})\simeq \langle\,\boldsymbol\mu_{r^\infty},\sqrt[r^\infty]{q_{E,v}}\,\rangle/\,q_{E,v}^\Z\ .\]
For any $v\in\mathcal{M}_F$ we choose a minimal Weierstrass equation for $E$.
Let $E_v$ be the reduction of $E$ modulo
$v$ and for any point $P\in E$ let $P_v$ be its image
in $E_v\,$.\\
For all basic facts about elliptic curves, the reader is referred
to Silverman's books \cite{Si1} and \cite{Si2}.

We remark that by increasing $q$ (if necessary) we can (and will) assume that $E$ is defined
over the field $F_0$ described in section \ref{FieldsNotations}.

\subsubsection{Duals} For $X$ a topological abelian group, we denote
its Pontrjagin dual by $X^{\vee}:=Hom_{cont}(X,\mathbb{C}^{\ast})$.
In the cases considered in this paper, $X$ will be a (mostly discrete)
topological $\Z_r$-module for some prime $r$, so that $X^{\vee}=Hom_{cont}(X,\Q_r/\Z_r)$ and it has a
natural structure of $\Z_r$-module.\\
The reader is reminded that to say that an $R$-module $X$ ($R$ any
ring) is cofinitely generated means that $X^{\vee}$ is a finitely
generated $R$-module. Since $(X^{\vee})^{\vee}\simeq X$, a module $X$ is $\Z_r$-cofinitely generated if and
only if it is the
direct sum of a finite ($r$-primary) abelian group with $(\Q_r/\Z_r)^t$ for
some $t\in\N$; in particular, letting $X_{div}$ be the divisible part of $X$, we see that $X/X_{div}$ is finite.

\subsection{Selmer groups}\label{Selflat} We shall deal with  torsion subschemes of the elliptic curve $E$.
Since $char\,F=p$, in order to deal with the $p$-torsion we need to consider flat cohomology of group schemes
to define the Selmer groups in that case.\\
For the basic theory of sites and cohomology on a site see \cite[Chapters II, III]{Mi1}.
We define our Selmer groups via flat cohomology (for the relation with classical
Galois cohomology see Remark \ref{flatetale} below) so, when we write a scheme $X$, we always mean the site $X_{fl}\,$.

Let $L$ be an algebraic extension of $F$ and $X_L:=Spec\,L$. For any positive integer $m$ the
group schemes $E[m]$ and $E$ define sheaves on $X_L$ (see \cite[II.1.7]{Mi1}): for example
$E[m](X_L):=E[m](L)$. Consider the exact sequence
\[ E[m]\iri E {\buildrel m\over\longrightarrow\!\!\!\!\!\rightarrow} E \]
and take flat cohomology to get
\[ E(L)/mE(L)\iri H_{fl}^1(X_L,E[m])\ri H_{fl}^1(X_L,E)\ . \]
In particular let $m$ run through the powers $r^n$ of a prime $r$. Taking direct limits
one gets an injective map (a ``Kummer homomorphism'')
\[ \k : E(L)\otimes \Q_r/\Z_r \iri \dl{n} H_{fl}^1(X_L,E[r^n])=: H_{fl}^1(X_L,E[r^\infty])\ . \]
As above one can build local Kummer maps for any place
$v\in\mathcal{M}_L$
\[ \k_v : E(L_v)\otimes \Q_r/\Z_r \iri H_{fl}^1(X_{L_v},E[r^\infty]) \]
where $X_{L_v}:=Spec\,L_v\,$.

\begin{defin}\label{Sel}
The \emph{$r$-part of the Selmer group} of $E$ over $L$, denoted by
$Sel_E(L)_r\,$, is defined to be
\[ Sel_E(L)_r:=Ker\left\{ H_{fl}^1(X_L,E[r^\infty])\ri\prod_{v\in\mathcal{M}_L}
H_{fl}^1(X_{L_v},E[r^\infty])/Im\,\k_v \,\right\} \]
where the map is the product of the natural restrictions between cohomology groups.\\
\end{defin}

The reader is reminded that if $L/F$ is finite then $Sel_E(L)_r$ is a cofinitely generated $\Z_r$-module.
Moreover the \emph{Tate-Shafarevich group}
$\ts(E/L)$ fits into the exact sequence
\[ E(L)\otimes \Q_r/\Z_r \iri Sel_E(L)_r \sri \ts(E/L)[r^\infty]\ .\]
According to the function field version of the Birch and
Swinnerton-Dyer conjecture, $\ts(E/L)$ is finite for any global function field
$L$. Applying to this last
sequence the exact functor $Hom(\cdot,\Q_r/\Z_r)$, it follows that
\[ rank_{\Z_r}\,Sel_E(L)_r^{\vee}=rank_\Z \, E(L) \]
(recall that cohomology groups, hence the Selmer groups, are endowed with
the discrete topology).

Fix a $\Z_l^d$-extension $\F/F$ with $l$ a prime different from $p$.
We will study the behaviour of the $r$-Selmer groups while $L$ varies
through the subextensions $F_n$ of $\F/F$. Such groups admit natural actions of
$\Z_r\,$, because of the torsion of $E$, and of $\G=Gal(\F/F)$.
Hence they are modules over the Iwasawa algebra $\Z_r[[\G]]$.
When $r=l$ this algebra is (noncanonically) isomorphic to the ring of
formal power series $\Zl[[T_1,\dots ,T_d]]$ (while, for $r\neq l$, $\Z_r[[\G]]$ is
more mysterious and we know virtually nothing about its structure).\\
In particular we will be concerned with the natural maps between
$\Z_r[[\G]]$-modules
\[ Sel_E(F_n)_r \ri Sel_E(\F)_r^{\G_n} \ .\]

\begin{rem}\label{flatetale}
{\em To define $Sel_E(L)_r$ (with $r\neq p$) we can also use the sequence
\[ E[r^n](\ov{F})\iri E(F^{sep}) {\buildrel r^n\over\longrightarrow\!\!\!\!\rightarrow} E(F^{sep}) \]
and classical Galois (=\'etale) cohomology since, in this case,
\[ H_{fl}^1(X_L,E[r^n])\simeq H_{et}^1(X_L,E[r^n])\simeq H^1(G_L,E[r^n](\ov{F})) \]
(see \cite[III.3.9]{Mi1}). To ease notations in this case
we shall write $H^i(L,\cdot)$ instead of
$H^i(G_L,\cdot)\simeq H^i_{fl}(X_L,\cdot)$ and write $E[n]$ for $E[n](\ov{F})$, putting
$E[r^\infty]:=\bigcup E[r^n]$. In this case the Kummer map
\[ \k : E(L)\otimes \Q_r/\Z_r \iri H^1(L,E[r^\infty]) \]
has an explicit description as follows. Let $\alpha\in E(L)\otimes \Qr/\Zr$ be represented by $\alpha= P\otimes
\frac{a}{r^k}\,$ ($a\in\Z$) and let $Q\in E(L^{sep})$ be such that $aP=r^kQ$. Then $\k(\alpha)=\varphi_\alpha\,$,
where $\varphi_\alpha(\sigma):=\sigma(Q)-Q$ for any $\sigma\in G_L\,$. }
\end{rem}

\section{Auxiliary lemmas}
We gather here the results which are needed for the proofs of the main theorems.
We start by giving a more precise description of $Im\,\k_v$ (following the path
traced by Greenberg in \cite{Gr1} and \cite{Gr2}). In our situation the local conditions
for the Selmer groups are easily seen to be often trivial (i.e., $Im\,\k_v=0$ in general),
a fact which is essentially due to $r\neq char\,F$.

\begin{prop}\label{Imkvl}
Let $L$ be the completion of an algebraic extension of $F_v$ and $r$
a prime different from $p$: then $E(L)\otimes \Qr/\Zr =0$ (i.e., the
Kummer map has trivial image).
\end{prop}

\begin{proof} This is an easy exercise: see e.g. \cite[Proposition 3.3]{BL}.\end{proof}

The following two lemmas deal with torsion points in abelian extensions of function fields
of characteristic $p$ both in the global and local case.

\begin{lemma}\label{GlobalTor}
Let $\F/F$ be a $\Z_l^d$-extension of function fields of characteristic $p>0$
and let $E/F$ be a non-isotrivial elliptic curve. Then the group $E(\F)_{tor}$ is finite.
\end{lemma}

\begin{proof}[Proof (sketch)] One proves a stronger statement: namely, that $E(L)_{tors}$
is finite for any abelian extension $L/F$.
Finiteness of $E[p^\infty](L)$ follows from the fact that points in $E[p^\infty]$ are inseparable over $F$
(a proof can be found e.g. in \cite[Proposition 3.8]{BLV}).
For the prime-to-$p$ part, it is shown in \cite[Theorem 4.2]{BLV} that the claim is a consequence of the following facts:\begin{itemize}
\item[{\bf 1.}] $Gal(F(E[r])/F)$ contains $SL_2(\bF_r)$ for almost all primes $r$;
\item[{\bf 2.}] $Gal(F(E[r^\infty])/F)$ contains $S_n$ for some $n$ (for any prime $r\neq p$)
where $S_n$ is the kernel of the natural reduction map $SL_2(\Z_r)\ri SL_2(\Z/r^n\Z)$.
\end{itemize}
Both statements follow from a theorem of Igusa \cite{Ig}. For a clear statement we refer to \cite{BLV},
where however appears the hypothesis that $F$ is global. So here we just show how to deduce {\bf 1} and
{\bf 2} in the case $F$ is not global.\\
Let $F_0$ be the global field described in section \ref{FieldsNotations}
and let $F'=F\cap F_0(E[r^\infty])$ (see the diagram below). The group $Gal(F'/F_0)$ is abelian because it is
a quotient of $Gal(F/F_0)\simeq Gal(\bF/\bF_q)$. Since $Gal(F'/F_0)\simeq
Gal(F_0(E[r^\infty])/F_0)/Gal(F_0(E[r^\infty])/F')$, one has that
$Gal(F_0(E[r^\infty])/F')$ contains the commutators of $Gal(F_0(E[r^\infty])/F_0)$.
By Igusa's theorem $Gal(F_0(E[r^\infty])/F_0)\supset S_n$ therefore
\[ S_{2n+2}\subset [S_n:S_n] \subset Gal(F_0(E[r^\infty])/F') \]
(for the inclusion on the right see e.g. \cite[Lemma 4.1]{BLV}).
Since $FF_0(E[r^\infty])=F(E[r^\infty])$ and the extensions $F/F'$ and $F_0(E[r^\infty])/F'$ are disjoint,
one gets $Gal(F(E[r^\infty])/F)\simeq Gal(F_0(E[r^\infty])/F')$ so
$Gal(F(E[r^\infty])/F)\supset S_{2n+2}$ as well.
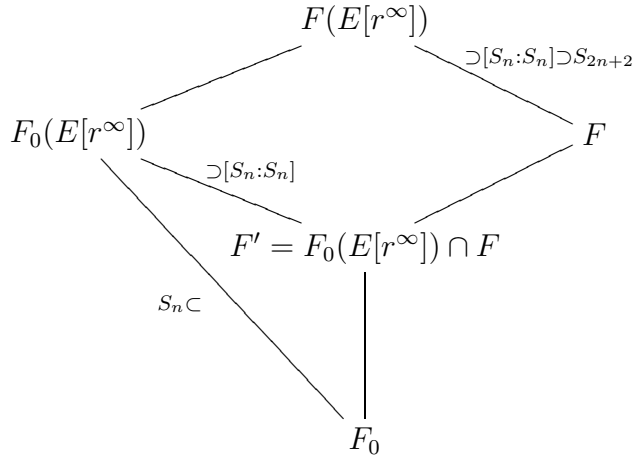
\[ \xymatrix{
\ & F(E[r^\infty])\ar@{-}[dl]\ar@{-}^{\ \ \ \supset[S_n:S_n]\supset S_{2n+2}}[dr] & \ \\
F_0(E[r^\infty])\ar@{-}^{\ \supset[S_n:S_n]}[dr]\ar@{-}_{S_n\subset\ }[dddr]  & \ & F \ar@{-}[dl] \\
\ & F'=F_0(E[r^\infty])\cap F\ar@{-}[dd] & \ \\
\ & \ & \\
\ & F_0 & \ } \]
This proves {\bf 2}. The same proof works for {\bf 1} as well (with $r$ in place of $r^\infty$),
remembering that $SL_2(\bF_r)$ is its own commutator subgroup for all primes $p\ge 5$.\end{proof}

\begin{lemma}\label{LocalTor}
Let $K$ be a field of characteristic $p$ complete with respect to
a discrete valuation $v$ and with residue field $\mathbb{K}\subset
\ov{\bF}_p\,$. Let $r$ be a prime different from $p$ and assume
that $\mathbb{K}$ does not contain $\bF_p^{(r)}$ (the
$\Z_r$-extension of $\bF_p\,$). Let $E/K$ be a non-isotrivial
elliptic curve. Then $E[r^\infty](K)$ is finite.
\end{lemma}

\begin{proof} Let $t$ be a uniformizer: then $K=\mathbb{K}((t))$ and exists $s$ such that
$E$ is defined over $K_0:=\bF_s((t))$. Since $K_0$ is a local field it is easy to see that
$E[r^\infty](K_0)$ is finite. Moreover since $\bF_p^{(r)}\not\subset\mathbb{K}$,
the Galois group $Gal(K/K_0)\simeq Gal(\mathbb{K}/\bF_s)$ contains no copies of $\Z_r\,$.\\
If $E[r^\infty](K)$ is infinite then choose an infinite sequence
of points $P_n\in E[r^n](K)$ such that $rP_{n+1}=P_n$ for any $n$.
Let $K'=K_0(\{P_n\}_{n\in \N})$ and $\mathcal P$ the subgroup of
$E[r^{\infty}]$ generated by the $P_n$'s. Then $K'/K_0$ is an
infinite extension and, since $K'\subset K$, one has
\[ Gal(K/K_0)\sri Gal(K'/K_0)\iri Aut(\mathcal P)\simeq \Z_r^*\ :\]
contradiction.\end{proof}

\begin{lemma}\label{H12Zld}
Let $\G\simeq \mathbb{Z}_l^d$ and $B$ a cofinitely generated discrete $\Zl$-module with a
continuous $\G$-action. Assume that there exists
a set $\g_1,\dots\g_d$ of independent topological generators of $\G$ such that $B^{\ov{\langle\g_1\rangle}}$ is finite.
Then, with $b:=\max\{\, |B/B_{div}|,|B^{\ov{\langle\g_1\rangle}}|\,\}$, one has
\[ |H^1(\G,B)|\le b^d\ and\ |H^2(\G,B)|\le b^{\frac{d(d-1)}{2}}\,. \]
\end{lemma}

\begin{proof} If $B$ is finite then $b=|B|$ and the proof
is in \cite[Lemma 4.1]{BL}. For the other case fix a set of
independent topological generators of $\G$ as above and put $\g:=\g_1\,$. Consider the exact sequence
\[ 0=B_{div}^{\ov{\langle\g\rangle}}\iri B_{div}\buildrel{\g-1}\over{-\!\!\!\!\longrightarrow} B_{div}
\sri B_{div}/(\g-1)B_{div} \]
(because of the hypothesis on $B$).
Taking duals one finds a sequence
\[ (B_{div}/(\g-1)B_{div})^\vee\iri (B_{div})^\vee\sri (B_{div})^\vee\simeq \Z_l^t \]
(for some finite $t$) and, counting ranks,
\[ {\rm rank}_\Zl (B_{div}/(\g-1)B_{div})^\vee=0 \ .\]
Therefore $(B_{div}/(\g-1)B_{div})^\vee$ is finite and, since $\Z_l^t$ has no
nontrivial finite subgroup, one finds
\[ B_{div}/(\g-1)B_{div}=0 \ .\]
Hence $B_{div}=(\g-1)B_{div}\subset (\g-1)B\subset B$ yields
\[ |B/(\g-1)B|\le |B/B_{div}| \ .\]
Now we use induction on $d$. For $d=1$ the equality $\G=\ov{\langle\g\rangle}$
implies $H^1(\G,B)\simeq B/(\g-1)B$ and $H^2(\G,B)=0$ (because $\Zl$
has $l$-cohomological dimension 1,
see \cite[Proposition 3.5.9]{NSW}).\\
For $d>1$ let $\G/\ov{\langle\g\rangle}=:\G'\simeq\Z_l^{d-1}$. The inflation
restriction sequence
\[ H^1(\G',B^{\ov{\langle\g\rangle}})\iri H^1(\G,B)\rightarrow H^1(\ov{\langle\g\rangle},B) \]
yields
\[ |H^1(\G,B)|\le |H^1(\G',B^{\ov{\langle\g\rangle}})|\ |H^1(\ov{\langle\g\rangle},B)|\le b^{d-1}b \ .\]
Moreover since $H^n(\ov{\langle\g\rangle},B)=0$ for any $n\ge 2$, the Hochschild-Serre spectral sequence
(see \cite[Theorem 2.1.5 and Exercise 5 page 96]{NSW}) gives an exact sequence
\[ H^2(\G',B^{\ov{\langle\g\rangle}})\rightarrow H^2(\G,B)\rightarrow H^1(\G',H^1(\ov{\langle\g\rangle},B)) \ .\]
By induction and the bound on $|H^1(\ov{\langle\g\rangle},B)|$ one has
\[ |H^2(\G,B)|\le | H^2(\G',B^{\ov{\langle\g\rangle}})|\ |H^1(\G',H^1(\ov{\langle\g\rangle},B))|\le \]
\[ \le b^{\frac{(d-1)(d-2)}{2}}b^{d-1} = b^{\frac{d(d-1)}{2}}\ .\]
\end{proof}

\begin{rem}\label{Remd=1}
\emph{Notice that if $d=1$ we have proved a slightly stronger
statement, namely that
\[ B^\G\ {\rm finite} \Longrightarrow |H^1(\G,B)|\le |B/B_{div}|\ .\]}
\end{rem}

To conclude we mention the version of Nakayama's Lemma we are going to use in what follows:
its proof (and further generalizations) can be found in \cite{BH}.

\begin{thm}\label{NakBH}
Let $\L$ be a compact topological ring with 1 and let $I$ be an ideal such that $I^n\rightarrow 0$.
Assume that $X$ is a profinite $\L$-module. If $X/IX$ is a finitely generated $\L/I$-module then
$X$ is a finitely generated $\L$-module and the number of generators of $X$ over $\L$
is at most the number of generators of $X/IX$ over $\L/I$.
Moreover if $\L=\Zl[[\G]]$, $I:=Ker\{\L\rightarrow \Z_l\,\}$ is the augmentation ideal
and $X/IX$ is finite then $X$ is $\L$-torsion.
\end{thm}

\section{Control theorems for $Sel_E(\F)_r$ ($r\neq p$)}\label{SecControll}
Before going on with the main theorems we describe the extensions we are going to deal with.
We recall that $\bF_p^{(r)}$ denotes the unique $\Z_r$-extension of $\bF_p\,$.

\begin{lemma}\label{ConstField}
For any prime $r\neq p$, the following statements are equivalent:
\begin{itemize}
\item[\bf{1.}] $\bF_p^{(r)}\subseteq \bF$;
\item[\bf{2.}] $\boldsymbol\mu_{r^\infty} \subset \bF(\boldsymbol\mu_r)$;
\item[\bf{3.}] $\Z_r \hookrightarrow Gal(\bF/\bF_p)$.
\end{itemize}
\end{lemma}

\begin{proof} Obvious, just recall that
\[ Gal(\ov{\bF}_p/\bF_p)\simeq \wh\Z :=\prod_r \Z_r \]
and
\[ \bF_p^{(r)} = \bF_p(\boldsymbol\mu_{r^\infty})^{Gal(\bF_p(\boldsymbol\mu_r)/\bF_p)}
\ .\]
\end{proof}

\begin{lemma}\label{LocExt}
Let $v$ be any place of $F$, $w$ a place of $\F$ dividing $v$ and $\G_v:=Gal(\F_w/F_v)$. One has that:
\begin{itemize}
\item[{\bf 1.}] if $\boldsymbol\mu_{l^{\infty}}\not\subset F_v$,
then
\[ \G_v \simeq \left\{ \begin{array}{cl} \Zl & {\rm if}\ v\ {\rm is\ inert} \\
0 & {\rm otherwise}\end{array}\right.\ ;\]
\item[{\bf 2.}] if
$\boldsymbol\mu_{l^{\infty}}\subset F_v$, then
\[ \G_v \simeq \left\{ \begin{array}{cl} \Zl & {\rm if}\ v\ {\rm is\ totally\ ramified} \\
0 & {\rm otherwise}\end{array}\right.\ .\]
\end{itemize}
\end{lemma}

\begin{proof} For any finite subextension $L/F_v$ of $\F_w/F_v$ we
have an exact sequence
$$I(L/F_v)\hookrightarrow Gal(L/F_v)\twoheadrightarrow Gal(\mathbb L/\bF_v)$$
where $I$ denotes the inertia subgroup. Since $\F_w/F_v$ is tamely
ramified, there is an injective homomorphism
$I(L/F_v)\hookrightarrow\bF_v^*$ (see e.g. \cite[IV, 2, Corollary
1 of Proposition 7]{Se}), hence
$|I(L/F_v)|\leq|\boldsymbol\mu_{l^{\infty}}(F_v)|$. There are two
cases.\\
{\bf Case 1}: $\boldsymbol\mu_{l^{\infty}}\not\subset F_v$. Since
$I(\F_w/F)$ is a submodule of the free $\Z_l$-module $\G_v$, it
follows from the boundedness of
$|\boldsymbol\mu_{l^{\infty}}(F_v)|$ and the equality
$I(\F_w/F_v)=\liminv I(L/F_v)$ that all these groups are trivial.
Therefore, either $\G_v\simeq Gal(\bF_v^{(l)}/\bF_v)$ and $\F_w$
is the constant field extension $\bF_v^{(l)}F_v$ or
$\F_w=F_v$.\\
{\bf Case 2}: $\boldsymbol\mu_{l^{\infty}}\subset F_v$. In this case $\bF_p^{(l)}\subset\bF_v$ and
$\bF_v$ has no $l$-extensions: hence either $\F_w=F_v$ or $\F_w/F_v$ is totally ramified. One
can apply Kummer theory to the classification of
$\Z_l$-extensions, as described in the Appendix. Let $t$ be a
uniformizer of the complete discrete valuation field $F_v$: from
$F_v^*=\bF_v^*\times t^{\Z}\times (1\text{-units})$ it follows
that the $l$-adic completion of $F_v^*$ is $t^{\Z_l}$, hence the
only $\Z_l$-extension is $F_v(\sqrt[l^{\infty}]{t})$.\end{proof}

\begin{prop}\label{ZlExt}
If $\bF_p^{(l)}\not\subset \bF$ then $F$ has a unique $\Z_l$-extension, namely
the constant field extension $\bF_p^{(l)}F$.
\end{prop}

For the proof, we remind the reader that $F$ is
the function field of a smooth, projective connected curve $\mathcal{C}$
defined over $\bF\,$. Remembering that $F=\bF F_0\,$, one sees that $\mathcal C$ can be obtained by base change
from a curve $\mathcal C_0$ defined over $\bF_q\,$.
Let $g$ be the genus of $\mathcal C_0$ and $\mathcal C\,$.

\begin{proof} Fix a geometric point $P$ of $\mathcal C$. By Lemma \ref{LocExt} one sees that a
$\Z_l^d$-extension $\F/F$ is everywhere unramified: therefore there is a
surjective morphism $\phi$ from the fundamental group $\pi_1(\mathcal C,P)$ to $Gal(\F/F)$.\\
We can assume that the point $P$ lies in $\mathcal C(\bF)$ (otherwise just take a finite
extension of $F$ whose constant field obviously still does not contain $\bF_p^{(l)}\,$).
Then we have a split exact sequence of fundamental groups
\[ \xymatrix{ \pi_1(\mathcal C\times \ov{\bF}_p, P) \ar@{^(->}[r] & \pi_1(\mathcal C,P) \ar@{->>}[r] &
G_\bF \ar@/_0,5cm/@{->}[l]} \ ,\]
that is, $\pi_1(\mathcal C,P)\simeq\pi_1(\mathcal C\times \ov{\bF}_p, P)\rtimes G_\bF$.
Since $Gal(\F/F)$ is abelian, the morphism $\phi$ factors through
$\pi_1(\mathcal C\times \ov{\bF}_p, P)^{ab}\rtimes G_\bF$ (notice that this semidirect product is a quotient
of $\pi_1(\mathcal C,P)\,$, since the $G_{\bF}$ action on $\pi_1(\mathcal C\times\ov{\bF}_p,P)$ preserves
the commutator subgroup).
It is well-known (see e.g. \cite[Proposition 9.1]{Mi3} together with \cite[XI, Th\'eor\`eme 2.1]{SGA1})
that one can identify $\pi_1(\mathcal C\times \ov{\bF}_p, P)^{ab}$ with the (full) Tate module of $Jac(\mathcal C)$.\\
Since $Gal(\F/F)$ is a pro-$l$ group (and the [pro]-primary-decomposition of a [profinite] abelian group
is preserved by automorphisms) the morphism $\phi$ factors further through
$T_l(Jac(\mathcal{C})) \rtimes G_\bF$. The following lemma shows that the maximal abelian quotient of
$T_l(Jac(\mathcal{C})) \rtimes G_\bF$ has the form $A\times G_\bF\,$, where $A$ is a finite group:
the proposition is an immediate consequence.\end{proof}

\begin{lemma}\label{ZlExtLemma}
If $\bF_p^{(l)}\not\subset \bF$ then the commutator subgroup of
$T_l(Jac(\mathcal{C})) \rtimes G_\bF$ has finite index in $T_l(Jac(\mathcal{C}))$.
\end{lemma}

\begin{proof} Since $G_\bF$  is abelian the commutators are contained in $T_l(Jac(\mathcal{C}))$.
To ease notation, shorten $T_l(Jac(\mathcal{C}))$ to $T$. We write the group law in $T\rtimes G_\bF$ as
$$(a,g)(b,h)=(a+gb,gh)$$
and let $\rho\colon G_\bF\rightarrow{\rm Aut}_\Zl(T)$ be the homomorphism corresponding to the action of $G_\bF$ on $T$. Then
$$(a,e)(0,h)(a,e)^{-1}(0,h)^{-1}=(a,h)(-a,e)(0,h^{-1})=(a-ha,h)(0,h^{-1})=(a-ha,e)$$
shows that to prove our claim it is enough to find $h\in G_\bF$ such that $(1-\rho(h))T$ has finite index in $T$.
Observe that since $T\simeq\Z_l^{2g}$ the operator $1-\rho(h)$ belongs to ${\rm End}_\Zl(T)\simeq M_{2g}(\Zl)$;
an easy reasoning shows that
$$[T:(1-\rho(h))T]=|\det(1-\rho(h))|_l^{-1}$$
(where $|\cdot|_l$ is normalized so that $|l|_l:=l^{-1}$). Hence we just need $\det(1-\rho(h))\neq0\,$.\\
Let $G_{\bF_q}^{(l)}$ and $G_{\bF}^{(l)}$ be respectively the maximal pro-$l$ subgroup of $G_{\bF_q}$ and
$G_{\bF}$: the hypothesis $\bF_p^{(l)}\not\subset \bF$ implies $[G_{\bF_q}^{(l)}:G_{\bF}^{(l)}]<\infty\,$.
Since all prime-to-$l$ subgroups of ${\rm Aut}_\Zl(T)\simeq GL_{2g}(\Zl)$ are finite
so is the index $[\rho(G_{\bF_q}):\rho(G_{\bF_q}^{(l)})]$. Hence there exists $h\in G_\bF^{(l)}$ such
that $\rho(h)=\rho(Frob_q^n)$ for some $n$ (where $Frob_q$ is the ``canonical'' generator of $G_{\bF_q}$).\\
The proof is concluded remarking the well-known fact that
$$\det(1-\rho(Frob_q^n))=| Jac(\mathcal C_0)(\bF_{q^n})|$$
and the right hand-side is not $0$.\end{proof}

We are now ready to prove two versions of the control theorem appropriate for our setting.

\subsection{The case $r=l$ with $\bF_p^{(l)}\not\subset\bF$}

\begin{thm}\label{Controll1}
Assume $\bF_p^{(l)}\not\subset\bF$.
Then the natural maps
\[ Sel_E(F_n)_l\rightarrow Sel_E(\F)_l^{\G_n} \]
have finite kernels and cokernels both of bounded order.
\end{thm}

\begin{proof} To ease notations, for any field $L$ let $\mathcal{G}(L)$
be the image of $H^1(L,E[l^\infty])$ in the product
\[ \prod_{w\in\mathcal{M}_L} H^1(L_w,E[l^\infty])/Im\,\k_w =
\prod_{w\in\mathcal{M}_L} H^1(L_w,E[l^\infty]) \]
(by Proposition \ref{Imkvl}).
We have a commutative diagram with exact rows
\[ \xymatrix {
Sel_E(F_n)_l \ar[d]^{a_n} \ar@{^{(}->}[r] &
H^1(F_n,E[l^\infty]) \ar[d]^{b_n} \ar@{->>}[r] &
\mathcal{G}(F_n)\ar[d]^{c_n}\\
Sel_E(\F)_l^{\G_n} \ar@{^{(}->}[r] &
H^1(\F,E[l^\infty])^{\G_n} \ar[r] &
\mathcal{G}(\F) } \]
and we are interested in $Ker\,a_n$ and $Coker\,a_n\,$.\\
By the Hochschild-Serre spectral sequence one gets
\[ Ker\,b_n\simeq H^1(\G_n,E[l^\infty](\F)) \]
and
\[ Coker\,b_n\subseteq H^2(\G_n,E[l^\infty](\F))\ .\]
By Lemma \ref{GlobalTor} the group $E[l^\infty](\F)$ is finite and by Proposition \ref{ZlExt} $\G_n\simeq \Z_l\,$.
So Lemma \ref{H12Zld} immediately gives
\[ |Ker\,b_n|\le |E[l^\infty](\F)| \quad {\rm and}\quad Coker\,b_n =0 \ . \]
By the snake lemma, this is enough to show that $Ker\,a_n$ is finite and bounded independently of $n$.\\
For $Coker\,a_n$ we need some control on $Ker\,c_n$ as well.
Obviously $Ker\,c_n$ embeds in the kernel of the natural map
\[ d_n : \prod_{v_n\in\mathcal{M}_{F_n}} H^1(F_{v_n},E[l^\infty])\longrightarrow
\prod_{w\in\mathcal{M}_{\F}} H^1(\F_w,E[l^\infty])\ .\]
For any $w|v_n$ we have a map
\[ d_w: H^1(F_{v_n},E[l^\infty])\longrightarrow H^1(\F_w,E[l^\infty]) \]
and $w_1\,,w_2|v_n$ imply $Ker\,d_{w_1}= Ker\,d_{w_2}\,$. Letting
$d_{v_n}$ be the product of the $d_w$'s for all the $w$'s dividing $v_n\,$, we have
$Ker\,d_{v_n}=\bigcap_{w_i|v_n}Ker\,d_{w_i}=Ker\,d_w$ for any $w|v_n$ and
\[ Ker\,c_n\subseteq Ker\,d_n =\prod_{v_n\in\mathcal{M}_{F_n}} Ker\,d_{v_n}\ .\]
By the inflation restriction sequence $Ker\,d_w = H^1(\G_{v_n},E[l^\infty](\F_w))$
(where $\G_{v_n}:=Gal(\F_w/F_{v_n})$ is independent of $w$ since $\G$ is abelian).

\noindent As seen in Lemma \ref{LocExt} one finds $\G_{v_n} =0$ or $\Zl$ and the latter
is the only nontrivial case. Moreover $\F_w/F_v$ is unramified (by Lemma \ref{ZlExt}):
therefore $\F_w\subset F_{v_n}^{unr}\,$, the maximal unramified extension of $F_{v_n}\,$.

\subsubsection{Places of good reduction} Assume $v_n$ is of good reduction.
By the criterion of N\'eron-Ogg-Shafarevich the field $F_{v_n}(E[l^\infty])$ is contained in $F_{v_n}^{unr}\,$.
The pro-$l$-part of $Gal(F_{v_n}^{unr}/F_{v_n})\simeq Gal(\ov{\bF}_p/\bF_{v_n})$ is isomorphic to
$\Z_l$ because $\bF_p^{(l)}\not\subset \bF$ yields $\bF_p^{(l)}\not\subset \bF_{v_n}$
(which is a finite extension of $\bF$).
Let $\varphi_l$ be a topological generator of the $\Z_l$-part of the Galois group
$Gal(F_{v_n}^{unr}/F_{v_n})$.
Since $H:=Gal(F_{v_n}^{unr}/F_{v_n})/\overline{\langle\varphi_l\rangle}$ has no $l$-primary part and
$E[l^\infty]$ is $l$-primary, the cohomology groups $H^i(H,E[l^\infty]^{\ov{\langle\varphi_l\rangle}})$
are trivial for $i\ge 1$. The Hochschild-Serre
spectral sequence provides an isomorphism
\[ H^1(Gal(F_{v_n}^{unr}/F_{v_n}),E[l^\infty])\simeq H^1(\overline{\langle\varphi_l\rangle},E[l^\infty])^H \ . \]
Note that the constant field of $F_{v_n,l}:=(F_{v_n}^{unr})^{\overline{\langle\varphi_l\rangle}}$ does not contain
$\bF_p^{(l)}$ because there is no $\Z_l$-extension between $F_{v_n}$ and
$F_{v_n,l}\,$. Therefore by Lemma \ref{LocalTor}, $E[l^\infty]^{\overline{\langle\varphi_l\rangle}}=E[l^\infty](F_{v_n,l})$
is finite. By Remark \ref{Remd=1} and the fact that $E[l^\infty]$ is divisible one has
$H^1(\overline{\langle\varphi_l\rangle},E[l^\infty])=0$, so $H^1(Gal(F_{v_n}^{unr}/F_{v_n}),E[l^\infty])$
is trivial too.
Since $\F_w\subset F_{v_n}^{unr}$, the inflation map
\[ H^1(\G_{v_n},E[l^\infty](\F_w))\iri H^1(Gal(F_{v_n}^{unr}/F_{v_n}),E[l^\infty]) \]
shows that
\[ Ker\,d_w = H^1(\G_{v_n},E[l^\infty](\F_w)) =0 \]
as well.

\subsubsection{Places of bad reduction} Let $\mathcal{R}_{n,i}$ be the (finite) set of primes of $F_n$
which are of bad reduction for $E$ and inert in $\F/F_n\,$. We recall that $\G_{v_n}\simeq \Z_l$ only if $v_n$ is inert
(otherwise $\G_{v_n}=0$); moreover $E[l^\infty](\F_w)^{\G_{v_n}}=E[l^\infty](F_{v_n})$ is finite by Lemma \ref{LocalTor}.
For a prime in $\mathcal R_{n,i}$, using Remark \ref{Remd=1} one immediately finds
\[ |Ker\,d_w| = |H^1(\G_{v_n},E[l^\infty](\F_w))| \le |E[l^\infty](\F_w)/E[l^\infty](\F_w)_{div}|\ . \]
Note that such bound actually depends on $v_n$ and not on $w$ so, to ease notations, we choose one prime
$w|v_n$ and we define
\[ \varepsilon(v_n):= |E[l^\infty](\F_w)/E[l^\infty](\F_w)_{div}| \ . \]
Therefore
\[ |Ker\,c_n|\le |Ker\,d_n|\le \prod_{v_n\in \mathcal{R}_{n,i}}
\varepsilon(v_n) \]
is finite and bounded as well.\end{proof}

\begin{rem}\label{BoundRank}
\emph{Recall that we are assuming that $E$ is a Tate curve at any (inert) place $v_n$ of bad reduction, so
\[ E[l^\infty](\F_w)_{div}=\left\{ \begin{array}{cl} 0 & \textrm{if}\ \boldsymbol\mu_l \not\subset\bF_v \\
\boldsymbol\mu_{l^\infty} & \textrm{if}\ \boldsymbol\mu_l \subset\bF_v \end{array}\right.\ .\]
Besides the Tate period $q_{E,v}$ has an $l^n$th root in $\F_w$ if and only if the $l$-adic valuation of
$ord_v(q_{E,v})$ is at least $n$.
Hence $E[l^\infty](\F_w)/E[l^\infty](\F_w)_{div}$ is a cyclic group of order
\[ \varepsilon(v)\le \frac{1}{|ord_v(q_{E,v})|_l} \ \]
(where $|\cdot|_l$ is the normalized $l$-adic absolute value).
Moreover (as in Lemma \ref{H12Zld}) one has a surjection
\[ E[l^\infty](\F_w)/E[l^\infty](\F_w)_{div}\sri E[l^\infty](\F_w)/(\gamma_{v_n}-1)E[l^\infty](\F_w)\simeq
H^1(\G_{v_n},E[l^\infty](\F_w)) \]
(where $\gamma_{v_n}$ is a topological generator of $\G_{v_n}$) which shows that $Ker\,d_w$ is generated by one element.}
\end{rem}

\begin{rem}\label{BoundOrder}
\emph{The uniform bounds provided by the theorem basically depend on the number of torsion points
and the places of bad reduction. Explicitly, letting $\mathcal{R}_i$ be the set of (inert) primes of $F$ 
of bad reduction for $E$, we found
\[ |Ker\,a_n|\le |E[l^\infty](\F)| \]
and
\[ |Coker\,a_n|\le \prod_{v\in\mathcal{R}_i} \varepsilon(v) \ .\]
Also, observe that $|E[l^\infty](\F)|$ is bounded by the number of
torsion points in the maximal abelian extension: so one could find
bounds depending only on $F$ and $E$.}
\end{rem}

\subsection{The case $r=l$ with $\bF_p^{(l)}\subset\bF$}

Notice that in this case, thanks to Lemmas \ref{ConstField} and \ref{LocExt}, only those places
$v$ such that $\boldsymbol\mu_l\subset\bF_v$ can ramify in $\F/F$; all the rest are totally split
(since $\bF_p^{(l)}\subset\bF$ there is no possibility for a $\Z_l$-extension of the constant field
corresponding to an inert $\Z_l$-extension of $F_v\,$).

\begin{thm}\label{Controll2}
Assume that $\bF_p^{(l)}\subset\bF$ and that only a finite number of places of $F$ ramify in $\F$.
Then the natural maps
\[ Sel_E(F_n)_l\rightarrow Sel_E(\F)_l^{\G_n} \]
have finite and bounded kernels and cofinitely generated cokernels
(of bounded corank over $\Zl$ when $d=1$).
\end{thm}

\begin{proof} Exactly as in Theorem \ref{Controll1},
we have a commutative diagram with exact rows
\[ \xymatrix {
Sel_E(F_n)_l \ar[d]^{a_n} \ar@{^{(}->}[r] &
H^1(F_n,E[l^\infty]) \ar[d]^{b_n} \ar@{->>}[r] &
\mathcal{G}(F_n)\ar[d]^{c_n}\\
Sel_E(\F)_l^{\G_n} \ar@{^{(}->}[r] &
H^1(\F,E[l^\infty])^{\G_n} \ar[r] &
\mathcal{G}(\F) } \]
with
\[ Ker\,b_n\simeq H^1(\G_n,E[l^\infty](\F)) \quad
{\rm and}\quad Coker\,b_n\subseteq H^2(\G_n,E[l^\infty](\F))\ .\]
Again by Lemma \ref{GlobalTor} the group $E[l^\infty](\F)$ is finite. Hence Lemma \ref{H12Zld} yields
\[ |Ker\,a_n\,|\le |Ker\,b_n\,|\le |E[l^\infty](\F)|^d \]
and
\[ |Coker\,b_n\,|\le  |E[l^\infty](\F)|^{\frac{d(d-1)}{2}}\ .\]
As before, for $Coker\,a_n$ we need some control on $Ker\,c_n$ and one gets it
by looking at the $Ker\,d_w=H^1(\G_{v_n},E[l^\infty](\F_w))$ for any $w|v_n\,$.

\subsubsection{Places of good reduction} Assume $v_n|v$ of good reduction.
By Lemma \ref{LocExt} we get $\G_{v_n}\simeq \Zl$ only if $v_n$ is ramified (otherwise it is 0 and
$Ker\,d_w$ is trivial).
Note that by the criterion of N\'eron-Ogg-Shafarevich
\[ E[l^\infty](\F_w)= E[l^\infty](F_v) \ .\]
Hence for a ramified place $v_n$ one has (with $\G_{v_n}=\ov{\langle \g_{v_n}\rangle}$)
\[ H^1(\G_{v_n},E[l^\infty](\F_w))=E[l^\infty](F_v)/(\g_{v_n}-1)E[l^\infty](F_v)=
 E[l^\infty](F_v) \]
which obviously has $\Zl$-corank $\le 2$ (notice that it can be equal to 2:
for example when $\bF=\ov{\bF}_p\,$).

\subsubsection{Places of bad reduction} Let $v_n$ be one of the (finitely many) primes of
bad reduction for $E$, lying above $v$. Since $\G_{v_n}$ is $\Zl$ or 0 it is easy to see that for these
ramified places
\[ corank_{\Zl} H^1(\G_{v_n},E[l^\infty](\F_w)) \le 2 \]
but we can be a bit more precise.\\
Assume $v_n$ is ramified (otherwise $Ker\,d_w=0$): by the theory of the Tate curve $E[l^\infty]\simeq
\langle\,\boldsymbol\mu_{l^\infty}\,,\,\sqrt[l^\infty]{q_{E,v}}\,\rangle/\,q_{E,v}^\Z$
where $q_{E,v}\in F_v$ is the Tate period (note that since $\boldsymbol\mu_{l^\infty}\subset\bF_v$
the set $E[l^\infty](\F_w)^{\G_{v_n}}=E[l^\infty](F_{v_n})$ is infinite and we cannot immediately apply
Lemma \ref{H12Zld}).
Besides $E[l^\infty](\F_w)=E[l^\infty]$. Therefore
\[  H^1(\G_{v_n},E[l^\infty](\F_w))\simeq H^1(\G_{v_n},\boldsymbol\mu_{l^\infty})\times
H^1(\G_{v_n},\sqrt[l^\infty]{q_{E,v}})\simeq
\boldsymbol\mu_{l^\infty} \] because $\G_{v_n}$ acts trivially on
$\boldsymbol\mu_{l^\infty}$ and $\sqrt[l^\infty]{q_{E,v}}$ is
divisible and such that $(\sqrt[l^\infty]{q_{E,v}})^{\G_{v_n}}$ is
finite (use Remark \ref{Remd=1}).\\

Let's divide the set of places ramified in $\F/F_n$ into $\mathcal{R}_{n,g}$
(consisting of primes where $E$ has good reduction)
and $\mathcal{R}_{n,b}$ (primes of bad reduction for $E$). Then all the above computations
lead to the bound
\[ corank_\Zl Coker\,a_n \le 2|\mathcal{R}_{n,g}\,| + |\mathcal{R}_{n,b}\,|\ . \]
Note that, if $d>1$, the number of ramified places is unbounded so
the coranks are unbounded as well, while for $d=1$ any ramified
place of $F$ can split only a finite number of times in
$\F$.\end{proof}

\begin{cor}\label{Controll2Fin}
In the setting of Theorem \ref{Controll2} assume that:\begin{itemize}
\item[1.] the ramified places are of good reduction for $E$;
\item[2.] $E[l^\infty](F_v)$ is finite for any ramified place $v$.
\end{itemize}
Then the natural maps $Sel_E(F_n)_l\rightarrow Sel_E(\F)_l^{\G_n}$
have finite (and bounded) kernels and finite cokernels (of bounded orders if
$d=1$).
\end{cor}

\begin{proof}
Just observe that the hypotheses yield
\[ Ker\,d_w = \left\{ \begin{array}{ll} 0 & {\rm if}\ v_n\ {\rm is\ unramified} \\
E[l^\infty](F_v) & {\rm otherwise} \end{array} \right.\ . \]
So one has $|Coker\,a_n\,|\le |E[l^\infty](\F)|^{\frac{d(d-1)}{2}}
\prod_{v_n\in \mathcal{R}_{n,g}} |E[l^\infty](F_v)|$.
\end{proof}

\begin{rem}
\emph{1. The assumption that only finitely many places ramify in
$\F/F$ is strictly necessary: see Example \ref{apex} in the appendix.\\
\indent 2. Hypotesis 2 in Corollary \ref{Controll2Fin} is often
satisfied. In case of good reduction, by the criterion of
N\'eron-Ogg-Shafarevich, we have $E[l^\infty](F_{v_n})\simeq
E_{v_n}[l^\infty](\bF_{v_n})=E_{v_n}[l^\infty]^G$, where
$G:=Gal(\bF_{v_n}(E_{v_n}[l^\infty])/\bF_{v_n})$. 
Let $\bF_q$ be the field of definition of $E_{v_n}$ and put
$G_0:=Gal(\bF_{v_n}(E_{v_n}[l^\infty])/\bF_q)$: as a quotient of 
$G_{\bF_q}$, $G_0$ is topologically generated by the Frobenius $Frob_q$. We consider the embedding 
$G_0\hookrightarrow{\rm Aut}(E_{v_n}[l^\infty])\simeq GL_2(\Zl)$: it's easy to see that $g\in G_0$ fixes
a finite number of points iff it has not 1 as an eigenvalue.
Assume that $Gal(\bF_{v_n}/\bF_q)\simeq\Zl$, so that if
$G_0$ has a prime-to-$l$ part, it must be $G$: in particular $G\neq\{1\}$ if 
the order of $Frob_q$ in ${\rm Aut}(E_{v_n}[l])$ does not divide $l$. Suppose besides
that the reduced curve ${\rm End}(E_{v_n})$ is an order $\mathcal O$ 
in a quadratic imaginary field $K$: then $Frob_q$ lies in ${\rm End}(E_{v_n})-\Z$ 
and it has eigenvalues $\{x,x^\tau\}$, $\tau$ a generator of $Gal(K/\Q)$.\footnote{We 
are just asking that $E_{v_n}$ is not supersingular: see \cite[V.3]{Si1}.} 
It follows that any $g\in G_0$ has eigenvalues $\{y,y^\tau\}$ for some 
$y\in\ov{\langle x\rangle}\subset(\mathcal O\otimes\Zl)^*$: 
in particular, if $l$ is not split in $K$, $y=1$ implies that $g$ is the identity.}
\end{rem}

Let $B$ be a cofinitely generated discrete $\Z_l$-module with a
continuous $\G$ action and denote $h_i(B)$ the number of
generators of $H^i(\G,B)$ ($i=1,2$). The same induction argument
as in Lemma \ref{H12Zld} shows that if $b$ is the
number of generators of $B$ then
\[ h_1(B) \le db \quad {\rm and}\quad
h_2(B) \le \frac{d(d-1)}{2}b\ .\]
One immediately finds the following corollaries (with identical proofs, so we only provide the first one).

\begin{cor}\label{Controll1Cor1}
In the setting (and with the notations) of Theorem \ref{Controll1} (and the subsequent remarks)
$Sel_E(\F)_l^\vee$ is a finitely generated
$\L$-module and
\[ rank_\L Sel_E(\F)_l^\vee \le corank_\Zl Sel_E(F)_l +|\mathcal{R}_i\,| \ .\]
Moreover if $Sel_E(F)_l$ is finite then $Sel_E(\F)_l^\vee$ is $\L$-torsion.
\end{cor}

This answers the analog of Question 1 and (some cases of) 2 in \cite{OT}.

\begin{cor}\label{Controll2Cor1}
In the setting (and with the notations) of Theorem \ref{Controll2}
$Sel_E(\F)_l^\vee$ is a finitely generated $\L$-module. Moreover
\[ rank_\L Sel_E(\F)_l^\vee \le corank_\Zl Sel_E(F)_l +2|\mathcal{R}_g\,| +
|\mathcal{R}_b\,|+h_2(E[l^\infty](\F)) \ ,\]
where $\mathcal{R}_g$ (resp. $\mathcal{R}_b\,$) is the set of ramified places
of $F$ of good (resp. bad) reduction for $E$ and, obviously, $h_2(E[l^\infty](\F))\le d(d-1)$.
\end{cor}

\begin{cor}\label{Controll2FinCor1}
In the setting of Corollary \ref{Controll2Fin}, if
$Sel_E(F)_l$ is finite then $Sel_E(\F)_l^\vee$ is a finitely generated torsion $\L$-module.
\end{cor}

\begin{proof} Let $\S$ be the Pontrjagin dual of
$Sel_E(\F)_l$ and let $I$ be the augmentation ideal of $\L$.
The quotient $\S/I\S$ is dual to $Sel_E(\F)_l^\G$ which is
cofinitely generated (resp. finite) by Theorem
\ref{Controll1} (resp. and the hypothesis on $Sel_E(F)_l\,$). Therefore Theorem \ref{NakBH} yields the
corollary. For the bound on the rank just use the exact sequences
\[ Sel_E(F)_l \ri Sel_E(\F)_l^\G \sri Coker\,a_0 \ ,\]
\[ Ker\,c_0\ri Coker\, a_0 \ri Coker\, b_0 = 0\]
and recall Remarks \ref{BoundRank} and \ref{BoundOrder}.
\end{proof}

\begin{rem}\label{RemEll}
\emph{For a computation of $rank_\L \S$ in the case $\bF=\ov{\bF}_p$ see
\cite[Propositions 2.5 and 3.4]{El}}
\end{rem}

\subsection{Applications} As well known, in case $d=1$ the structure of the dual of Selmer groups 
can be used to control the growth of Mordell-Weil ranks in the tower of extensions 
between $F$ and $\F$ and to formulate an ``Iwasawa Main Conjecture''.

\subsubsection{Mordell-Weil ranks}\label{MWR}
In \cite[Theorem 1.1]{Sh} Shioda proves that the group $E(F)$ is finitely generated for any function field $F$ with
algebraically closed constant field (of course this covers
the case of the $\Z_l$-extension $\bF_p^{(l)}F$ as well).
Our Corollary \ref{Controll2FinCor1} provides a new family of extensions for which $E(\F)$ is finitely generated.

\begin{cor}\label{FinGenE}
In the setting of Corollary \ref{Controll2Fin} assume that $\F/F$ is a $\Zl$-extension and that $Sel_E(F)_l$ is
finite. Then $E(\F)$ is finitely generated.
\end{cor}

\begin{proof} (More details can be found in \cite[Theorem 1.3 and Corollary 4.9]{Gr2})
Let $\S$ be the dual of $Sel_E(\F)_l\,$:
by Corollary \ref{Controll2FinCor1}, $\S$ is a finitely generated torsion $\L$-module.
By the well-known structure theorem for such modules there is a pseudo-isomorphism
\[ \S\sim \bigoplus_{i=1}^s \Zl[[T]]/(f_i^{e_i})\,. \]
Let $\lambda=\deg \prod f_i^{e_i}\,$: then ${\rm rank}_\Zl \S=\lambda$ and, taking duals, one gets
\[ (Sel_E(\F)_l)_{div}\simeq (\Ql/\Zl)^\lambda \ .\]
By Corollary \ref{Controll2Fin}, for any $n$, one has
\[ (Sel_E(F_n)_l)_{div}\simeq (\Ql/\Zl)^{t_n}\ {\rm with}\ t_n\le\lambda\ . \]
Hence
\[ (\Ql/\Zl)^{rank\,E(F_n)} \simeq E(F_n)\otimes \Ql/\Zl \iri (Sel_E(F_n)_l)_{div} \]
yields $rank\,E(F_n)\le t_n\le \lambda$ for any $n$, i.e. such ranks are bounded.\\
Choose $m$ such that $rank\,E(F_m)$ is maximal and let $t=|E(\F)_{tor}|$.
Using the fact that $E(\F)/E(F_m)$ is a torsion group one proves that $tP\in E(F_m)$ for all $P\in E(\F)$
and multiplication by $t$ gives a homomorphism $\varphi_t:\ E(\F)\rightarrow E(F_m)$ whose image is finitely
generated and whose kernel is the finite group $E(\F)_{tor}\,$.
Hence $E(\F)$ is indeed finitely generated.\end{proof}

\subsubsection{Iwasawa Main Conjecture}\label{IMC}
When $F$ is a global field (and, necessarily, $d=1$ and $\F=\bF_p^{(l)}F$),
our control theorem may be used, as classically, as a first step towards the algebraic side for a Main Conjecture.
As for the analytic side, the best candidate we know of has been provided by P\'al. In \cite{Pa}, he constructs an
element $\mathcal L_\infty(E)$ in the Iwasawa algebra $\Z[[G_{\infty}]]\otimes\Q$ (where $G_\infty$ is the Galois
group of the maximal abelian extension of $F$ unramified outside a fixed place where $E$ has split multiplicative
reduction). He is then able to prove an interpolation formula connecting $\mathcal L_\infty(E)$ to a special value
of the classical Hasse-Weil $L$-function of $E$ (\cite[Theorem 1.6]{Pa}). Now, since $\G$ is a quotient of $G_\infty$,
there is a natural map $\pi\colon\Z[[G_{\infty}]]\otimes\Q\rightarrow\Zl[[\G]]\otimes\Q\,.$
The element $\mathcal L_\G(E):=\pi(\mathcal L_\infty(E))$ would then be a natural candidate for a generator of the characteristic
ideal of $Sel_E(\F)^\vee_l$.

Support for such a conjecture comes from recent work of Ochiai and Trihan \cite{OT}. By means of techniques of syntomic
cohomology, they are able to prove an Iwasawa Main Conjecture for a semistable abelian variety $A/F$ and the $\Zp$-extension $F_\infty^{(p)}:=\bF_p^{(p)}F$ \cite[Theorem 1.4]{OT}. It is not known yet what is the relation (if any) between
Pal's $\mathcal L_\infty(E)$ and Ochiai-Trihan's $\mathcal L_{A/F_\infty^{(p)}}$ (but see \cite[Remark 3.2]{OT}).

We also remark that Ochiai-Trihan \cite{OT0} are able to prove that their Selmer dual is always torsion
(a necessary condition to have a non-zero characteristic ideal). So one expects the analog to be true for
our $Sel_E(\F)^\vee_l$ as well.

\subsection{The case $r\neq l,p$}
The $r$-part of Selmer groups behaves well in a $\Z_l^d$-extension: indeed it is easy to see that

\begin{thm}\label{Controlr}
The natural maps $Sel_E(F_n)_r\rightarrow Sel_E(\F)_r^{\G_n}$ are isomorphisms.
\end{thm}

\begin{proof} We use the same diagram of Theorem \ref{Controll1}, only changing
$l$-torsion with $r$-torsion points (since $r\neq p$ we can still use Galois cohomology).
The proof goes on in the same way noting that
\[ Ker\,b_n = H^1(\G_n,E[r^\infty](\F)) = 0\ ,\]
\[ Coker\,b_n \subseteq H^2(\G_n,E[r^\infty](\F)) = 0\ ,\]
\[ Ker\,d_w = H^1(\G_{v_n},E[r^\infty](\F_w)) = 0 \]
because $E[r^\infty](\F)$ and $E[r^\infty](\F_w)$ are $r$-primary while
$\G_n$ and $\G_{v_n}$ are pro-$l$-groups.\end{proof}

The consequences of this theorem on the structure of $Sel_E(\F)_r$ as a $\Z_r[[\G]]$-module will be given in the
next section together with the results on $Sel_E(\F)_p$ (see Corollary \ref{ControlpCor}).

\section{Control theorem for $Sel_E(\F)_p$}\label{SecControlp}
In this section we shall work with the $p$-torsion; so
we need flat cohomology, as explained in section \ref{Selflat}, and we shall
follow the notations given there. \\
As before, it is convenient to write $\F=\bigcup F_n$ with $F_n/F$ finite and
$F_n\subset F_{n+1}\,$.\\

\begin{thm}\label{Controlp}
The natural maps $Sel_E(F_n)_p\longrightarrow Sel_E(\F)_p^{\G_n}$
are isomorphisms.
\end{thm}

\begin{proof} We start by fixing the notations which will be used throughout the proof.\\ Let
$X_n:=Spec\,F_n\,$, $\X:=Spec\,\F\,$, $X_{v_n}:=Spec\,F_{v_n}$ and $\X_{w}:=Spec\,\F_{w}\,$.
To ease notations, let
\[ \mathcal{G}(X_n):=Im\,\left\{H_{fl}^1(X_n,E[p^\infty])\ri
\prod_{v_n\in\mathcal{M}_{F_n}}\,H_{fl}^1(X_{v_n},E[p^\infty])/Im\,\k_{v_n}\,\right\} \]
(analogous definition for $\mathcal{G}(\X)\,$).

\noindent Just like in the previous section we have a diagram
\[ \xymatrix{
Sel_E(F_n)_p \ar[d]^{a_n} \ar@{^{(}->}[r] & H_{fl}^1(X_n,E[p^\infty]) \ar[d]^{b_n} \ar@{->>}[r] &
\mathcal{G}(X_n)\ar[d]^{c_n}\\
Sel_E(\F)_p^{\G_n} \ar@{^{(}->}[r] & H_{fl}^1(\X,E[p^\infty])^{\G_n} \ar[r] & \mathcal{G}(\X) \ .} \]

\subsection{The map $b_n\,$.} The map $\X\ri X_n$ is a Galois covering with Galois group $\G_n\,$. In this context
the Hochschild-Serre spectral sequence holds by \cite[III.2.21 a),b) and III.1.17 d)]{Mi1}.
Therefore one has an exact sequence
\[   H^1(\G_n,E[p^\infty](\F))\iri  H_{fl}^1(X_n,E[p^\infty])\rightarrow
  H_{fl}^1(\X,E[p^\infty])^{\G_n} \rightarrow H^2(\G_n,E[p^\infty](\F))  \]
which fits in the diagram above (note that the first and last elements are Galois cohomology groups).\\
Since $E[p^\infty](\F)$ is a finite $p$-primary group (by Lemma \ref{GlobalTor}) and $\G_n$ is a pro-$l$-group, one has
\[ H^i(\G_n,E[p^\infty](\F))=0\quad (i=1,2) \]
and $Ker\,b_n=Coker\,b_n=0$ as well.

\subsection{The map $c_n\,$.} First of all we note that $Ker\,c_n$ embeds into the kernel of the map
\[  d_n : \prod_{v_n\in\mathcal{M}_{F_n}}\,H_{fl}^1(X_{v_n},E[p^\infty])/Im\,\k_{v_n}
\longrightarrow \prod_{w\in\mathcal{M}_{\F}}\, H_{fl}^1(\X_{w},E[p^\infty])/Im\,\k_{w} \]
and we only consider the maps
\[  d_w : H_{fl}^1(X_{v_n},E[p^\infty])/Im\,\k_{v_n}
\longrightarrow H_{fl}^1(\X_w,E[p^\infty])/Im\,\k_{w} \]
separately. Observe that:
\begin{itemize}
\item[{\bf 1.}] for any $v_n$ there are as many maps $d_w$ as many
primes $w$ of $\F$ dividing $v_n$ but all these maps have isomorphic kernels;
\item[{\bf 2.}] $Ker\,c_n\subseteq
\prod_{v_n\in\mathcal{M}_{F_n}} \bigcap_{w|v_n} Ker\,d_w\,$.
\end{itemize}

\noindent The Kummer exact sequence yields a diagram
\[ \xymatrix{
H_{fl}^1(X_{v_n},E[p^\infty])/Im\,\k_{v_n}\ar@{^{(}->}[r]\ar[d]^{d_w} &
H_{fl}^1(X_{v_n},E)[p^\infty]\ar[d]^{h_w} \\
H_{fl}^1(\X_w,E[p^\infty])/Im\,\k_w\ar@{^{(}->}[r] & H_{fl}^1(\X_w,E)[p^\infty] \ .} \]
Again $\X_w \ri X_{v_n}$ is a Galois covering so the Hochschild-Serre spectral sequence implies
\[ Ker\,d_w\iri Ker\,h_w\simeq H^1(\G_{v_n},E(\F_w))[p^\infty] = \dl{k} H^1(\G_{v_n},E(\F_w))[p^k] \ . \]
But $H^1(\G_{v_n},E(\F_w))[p^k]=0$ because it consists of the $p^k$-torsion of the cohomology of
a pro-$l$-group.\\
This yields $Ker\,c_n=0$ and therefore $a_n$ is an isomorphism.\end{proof}

\subsection{Structure of $Sel_E(\F)_r$ for $r\neq l$.} The Selmer groups $Sel_E(\F)_r$
are modules over the ring $\Z_r[[\G]]$ and,
to apply the generalized Na\-kaya\-ma's Lemma of \cite{BH} (i.e. Theorem \ref{NakBH} above), we need an ideal $J$ of $\Z_r[[\G]]$
such that $J^n\rightarrow 0$.
The classical augmentation ideal $I$ does not verify this condition since $I=I^2$
(see \cite[Lemma 3.7]{BL}).\\
Anyway we can use the ideal $rI$ to obtain a partial description of $Sel_E(\F)_r\,$.
We need the following (detailed proof in \cite[Lemma 3.8]{BL}).

\begin{lemma}\label{rIdual}
Let $M$ be a discrete $\Z_r[[\G]]$-module and $m_r: M\rightarrow M$ the multiplication by $r$. Then
\[ M^{\vee}/rIM^{\vee}\simeq (m_r^{-1}(M^\G))^\vee = (M^\G +M[r])^\vee \]
(where $M[r]$ is the $r$-torsion of $M$).
\end{lemma}

\begin{proof} Let $N=M^\vee$ so that $N$ is a $\Z_r[[\G]]$-module. Via the dual of the
natural projection map $\pi : N\sri N/rIN$ one sees that
\[ (N/rIN)^{\vee}\simeq m_r^{-1}((N^{\vee})^\G)\ ,\]
which yields
\[  M^{\vee}/rIM^{\vee}\simeq (m_r^{-1}(M^\G))^{\vee}\ .\]
Since $H^1(\G,M[r])=0$ one has $m_r(M)^\G=m_r(M^\G)$ and
can conclude noting that
\[ m_r^{-1}(M^\G) = m_r^{-1}(m_r(M^\G))=M^\G+M[r]\ .\]
\end{proof}

\begin{cor}\label{ControlpCor}
Assume that both $Sel_E(F)_r$ and $Sel_E(\F)_r[r]$ are finite. Then $Sel_E(\F)_r^{\vee}$
is a finitely generated $\Z_r[[\G]]$-module.
\end{cor}

\begin{proof} By the previous lemma with $M=Sel_E(\F)_r$ one has
\[ Sel_E(\F)_r^{\vee}/rI Sel_E(\F)_r^{\vee} \simeq
(Sel_E(\F)_r^\G+Sel_E(\F)_r[r])^{\vee}  \]
so this quotient is finite by hypothesis and Theorems \ref{Controlr} or \ref{Controlp}.
Then Theorem \ref{NakBH} yields our corollary.\end{proof}

In the corollary it would be enough to assume that $Sel_E(F)_r$ and $Sel_E(\F)_r[r]$ are
cofinitely generated modules over $\Z_r[[\G]]/rI \Z_r[[\G]]$. Unfortunately
even with the stronger assumption of finiteness we can't go further
(i.e., we are not able to see whether $Sel_E(\F)_r^{\vee}$ is
a torsion $\Z_r[[\G]]$-module or not) due to our lack of understanding of the structure of $\Z_r[[\G]]$-modules
even for simpler $\G$'s like for example $\G\simeq \Zl\,$.

\appendix
\section{$\Zl$-extensions of a field}

Let $F$ be a field, on which we assume only that
$\boldsymbol\mu_{l^{\infty}}\subset F$, with $l\neq char(F)$ a
prime. Everything is taking place in a fixed separable closure
$F^{sep}$. The goal is to describe the set of all
$\Z_l^d$-extensions of $F$ in $F^{sep}$.\\
\indent Define $\wh{F^*}$ as the $l$-adic completion of $F^*$: that is,
$\widehat{F^*}:=\liminv F^*/(F^*)^{l^n}$. This is a topological
$\Z_l$-module (each quotient $F^*/(F^*)^{l^n}$ is given the
discrete topology) and the natural map $F^*\rightarrow\wh{F^*}$
has dense image.\\
Let $V:=\Q_l\otimes_{\Zl}\wh{F^*}$. Then $V$ is a topological
$\Q_l$-vector space, complete and locally convex, with a
distinguished lattice $\wh{F^*}$ (more precisely, $V$ is a Banach
space over $\Q_l$, with the norm induced by taking $\wh{F^*}$ as
unit ball). The natural map $\wh{F^*}\rightarrow V$ is an
injection.\\

The reader is reminded that, if $W$ is a vector space, the
Grassmannian ${\rm Grass}_d(W)\subset\mathbb{P}(\Lambda^dW)$ is
the set of all $d$-dimensional subspaces of $W$.

\begin{thm} The set of $\Z_l^d$-extensions of $F$ is in bijection with ${\rm
Grass}_d(V)$.
\end{thm}

\begin{proof} By the assumption on $\boldsymbol\mu_{l^{\infty}}$, we have that
$\Z_l(1):=\liminv\boldsymbol\mu_{l^n}$ is isomorphic to $\Z_l$ as
$G_F$-module. Hence a $\Z_l^d$-extension $\F/F$ is uniquely
determined by the kernel of a continuous homomorphism
$G_F\rightarrow\Z_l(1)^d$ with image a rank $d$ submodule
($\Z_l(1)$ is given the profinite topology).\\
\indent We have
$${\rm Hom}_{cont}(G_F,\Z_l(1)^d)\simeq
{\rm Hom}_{cont}(G_F,\Z_l(1))^d\simeq\big(\liminv {\rm
Hom}(G_F,\boldsymbol\mu_{l^n})\big)^d\simeq\wh{F^*}^d$$ where all
isomorphisms\footnote{These are isomorphisms of topological
groups, giving to ${\rm Hom}_{cont}(G_F,\bullet)$ the compact open
topology. Notice that since $\boldsymbol\mu_{l^n}$ is discrete so
is also ${\rm Hom}(G_F,\boldsymbol\mu_{l^n})$.} are almost
tautological but the last one, which comes from Hilbert 90 and the
observation that the diagram
\[ \begin{CD} F^*/(F^*)^{l^{n+1}} @>>> {\rm Hom}(G_F,\boldsymbol\mu_{l^{n+1}}) \\
@VVV @VVV \\
F^*/(F^*)^{l^n} @>>> {\rm Hom}(G_F,\boldsymbol\mu_{l^n})
\end{CD} \]
commutes. Here, for any $n$, horizontal maps are the Kummer homomorphisms
sending $a\in F^*/(F^*)^{l^n}$ to
$\sigma\mapsto\frac{\sigma\sqrt[l^n]{a}}{\sqrt[l^n]{a}}$
and the right-hand vertical map is induced by raising-to-$l$:
$\boldsymbol\mu_{l^{n+1}}\rightarrow\boldsymbol\mu_{l^n}$.\\
That is, any continuous homomorphism $G_F\rightarrow\Zl(1)^d$ is of
the form $\langle\cdot,x\rangle=\lim\langle\cdot,x_n\rangle_n$ for
some $x=(x_{i,n})\in\wh{F^*}^d$, where
$\langle\cdot,\cdot\rangle_n \colon G_F\times
(F^*/(F^*)^{l^n})^d\rightarrow\boldsymbol\mu_{l^n}^d$ is the
$l^n$th level Kummer pairing,
$\langle\sigma,y\rangle_n:=(\frac{\sigma\sqrt[l^n]{y_1}}{\sqrt[l^n]{y_1}},
..., \frac{\sigma\sqrt[l^n]{y_d}}{\sqrt[l^n]{y_d}})$.\\
Let $\F_x\subset F^{sep}$ be the fixed field of
$\ker\langle\cdot,x\rangle$ and $B_x$ the closure of the subgroup
of $\wh{F^*}$ generated by $x_1,...,x_d$. It is well-known that
$F_{x,n}:=F(\sqrt[l^n]{x_{1,n}},...,\sqrt[l^n]{x_{d,n}})$ is the
fixed field of $\ker\langle\cdot,x_n\rangle_n$ and that
$Gal(F_{x,n}/F)\simeq G_F/\ker\langle\cdot,x_n\rangle_n$ is the
dual of $B_x/(\wh{F^*})^{l^n}$. It follows that
$\F_x=\bigcup_nF_{x,n}$ (since
$\ker\langle\cdot,x\rangle=\cap\ker\langle\cdot,x_n\rangle_n$) and
that $Gal(\F_x/F)$ is (non-canonically) isomorphic to
$B_x\simeq\liminv B_x/(\wh{F^*})^{l^n}$ (because any finite
abelian group is non-canonically isomorphic to its
dual).\\
In the same way, one sees that $\F_x=\F_y$ if and only if
$B_x\otimes\Q_l=B_y\otimes\Q_l$.

The theorem follows.
\end{proof}

\begin{ex} \label{apex} {\em Let
$F=\overline\bF_p(T)$ and choose a family $a_i\in\overline\bF_p$,
$i\in \N$ and $a_i\neq a_j$ if $i\neq j$. Put $\pi_i:=T+a_i$ and
consider the sequence
\[ x_1=\pi_1 \ ,\  x_2=x_1\pi_2^l\  ,\  x_3=x_2\pi_3^{l^2}\dots \]
\[ x_{n+1}=x_n\pi_{n+1}^{l^n}\ .\]
The elements $x_i$ provide a $\Z_l$-extension
\[ \F_x=\bigcup_{n\in\N} F(\sqrt[l^n]{x_n}) \]
ramified at all the $\pi_i$'s.}
\end{ex}

\noindent A. Bandini\\
Universit\`a della Calabria - Dipartimento di Matematica\\
via P. Bucci - Cubo 30B - 87036 Arcavacata di Rende (CS) - Italy\\
bandini@mat.unical.it

\noindent I. Longhi\\
National Taiwan University - Department of Mathematics\\
No. 1 section 4 Roosevelt Road - Taipei 106\\
Taiwan\\
longhi@math.ntu.edu.tw

\end{document}